\documentclass[11pt]{article}

\usepackage{amssymb}
\usepackage{a4wide}

\newtheorem{theorem}{Theorem}

\newcommand{\Gal}{\mathrm{Gal}}
\newcommand{\GL}{\mathrm{GL}}
\newcommand{\Pic}{\mathrm{Pic}}
\newcommand{\End}{\mathrm{End}}
\newcommand{\Aut}{\mathrm{Aut}}
\newcommand{\tors}{\mathrm{tors}}
\renewcommand{\deg}{\mathrm{deg}}

\newcommand{\ZZ}{\mathbb{Z}}
\newcommand{\QQ}{\mathbb{Q}}

\newcommand{\CC}{\mathbb{C}}

\newcommand{\qed}{\quad $\square$\bigskip}

\begin{document}

\title{Images of isogeny classes on modular elliptic curves}
\author{Florian Breuer\footnote{The author would like to thank the 
Max-Planck-Institut f\"ur Mathematik, Bonn, where this paper was written.}  \\
%
%
{\small
Department of Mathematics, University of Stellenbosch, Stellenbosch
7600, South Africa.} \\  
{\small fbreuer@sun.ac.za}}
\maketitle

\begin{abstract}
Let $K$ be a number field and $E/K$ a modular elliptic curve, with
modular parametrization $\pi : X_0(N) \longrightarrow E$ defined over $K$.  The
purpose of this note is to study the images in $E$ of classes of isogenous
points in $X_0(N)$.
\end{abstract}

Let $\pi : X_0(N) \rightarrow E$ be as above, and denote by $\bar{K}$ an
algebraic closure of $K$.

\begin{theorem}
Let $S \subset X_0(N)(\bar{K})$ be an infinite set of points corresponding to
elliptic curves which all  lie in one isogeny class, but which are not
isogenous to $E$ itself. Then the subgroup of $E(\bar{K})$ generated by
$\pi(S)$ has infinite rank and finite torsion.
\end{theorem}

\paragraph{Proof.}
Write $S=\{x_0,x_1,\ldots\}$ and $y_i := \pi(x_i)\in E(\bar{K})$ for $i\geq
0$.  We first show that $\langle \pi(S) \rangle$ is not finitely generated, and
then that it has finite torsion.

Suppose that $\langle \pi(S) \rangle$ is finitely generated. Then $\langle
\pi(S) \rangle \subset E(L)$ for some number field $L$, which we may extend to
include $K$. Now $G_L := \Gal(\bar{L}/L)$ acts on each fiber $\pi^{-1}(y_i)$,
from which follows that
\begin{equation}\label{eq1}
|G_L\cdot x_i| \leq \deg(\pi),\qquad\forall i\geq 0.
\end{equation}

Denote by $E_i$ the elliptic curve corresponding to $x_i$ for each $i\geq 0$.
It is isogenous to $E_0$. We now consider two cases. 

{\bf (i)} If $E_0$ has complex multiplication, then each $\End(E_i)$ is an 
order of conductor $f_i$ in a fixed quadratic imaginary field $F$. 
We denote by $h_F$ the class number of $F$. Then we have
\begin{eqnarray*}
|G_L\cdot x_i| & \geq & |\Pic(\End(E_i))|/2[L : \QQ] 
\quad\mbox{(by \cite[Chap 10, Theorem 5]{Lang})}\\
& \geq & \frac{h_F}{12[L:\QQ]}\cdot f_i\prod_{p|f_i}\left(1-\frac{1}{p}\right) 
\quad\mbox{(by \cite[Chap 8, Theorem 7]{Lang}),}
\end{eqnarray*}
which tends to $\infty$ as $i \rightarrow\infty$, thus contradicting
(\ref{eq1}).
 
{\bf (ii)} Now suppose that $E_0$ does not have complex multiplication. We may
write $E_i = E_0/C_i$, with $C_i \subset E_0$ a cyclic subgroup of order $n_i$.
Consider the Galois representations
\[
\rho_{n_i} : G_L \longrightarrow \Aut(E_0[n_i]) \cong \GL_2(\ZZ/n_i\ZZ)
\]
attached to $E_0$. From \cite[Th\'eor\`eme 3$^\prime$]{Serre} follows that 
there exists a constant
$d_0$, depending only on $E_0$ and on $L$, such that the image of $\rho_{n_i}$
has index at most $d_0$. Thus
\[
|G_L\cdot x_i| \geq
|\Aut(E_0[n_i])\cdot C_i|/d_0
= \psi(n_i)/d_0,
\]
where $\psi(n_i)=n_i\prod_{p|n_i}(1+1/p) \geq n_i$ is the number of cyclic
subgroups of order $n_i$ in $E_0[n_i]$. This again contradicts (\ref{eq1}), and
it follows that $\langle \pi(S) \rangle$ is not finitely generated. Notice 
that at this point we have not yet used the assumption that the $E_i$ are
not isogenous to $E$ itself.

We now show that $\langle \pi(S) \rangle$ has finite torsion. Let $K_0 \supset
K$ be a number field over which $E_0$ is defined, then every $E_i$ is defined
over $L_0 = K_0(E_{0,\tors})$. From the Weil pairing follows that
$K_0(\mu_{\infty})\subset L_0$. From \cite[Th\'eor\`eme
6$^{\prime\prime\prime}$]{Serre}  and \cite[Satz 4]{Faltings} follows that
$L_0\cap K_0(E_{\tors})$ is a finite extension of $K_0(\mu_{\infty})$, as $E$
and $E_0$ are not isogenous. Therefore we may write $L_0\cap K_0(E_{\tors}) =
L(\mu_{\infty})$ for some number field $L$. Now from \cite{Ribet} follows that
$E_{\tors}(L(\mu_{\infty}))$ is finite, yet  $\langle \pi(S) \rangle_{\tors}
\subset E_{\tors}(L(\mu_{\infty}))$, which completes our proof. 
\qed

What happens with images of points isogenous to $E$ itself? Here it is
conceivable that the image has infinite torsion, but the following result shows
that if $S$ contains an infinite chain of cyclic $m$-isogenies 
$x_1\stackrel{m}{\rightarrow} x_2\stackrel{m}{\rightarrow}\cdots$ for $m$ 
sufficiently large, then infinitely many of the $\pi(x_i)$'s must be points of
infinite order.

\begin{theorem}
Let $m \geq\max(2,\deg(\pi))$. Then there exist only finitely many pairs of
torsion points $y_1,y_2\in E_{\tors}(\CC)$ which possess preimages
$x_1\in\pi^{-1}(y_1)$ and $x_2\in\pi^{-1}(y_2)$ corresponding to elliptic
curves  $E_1$ and $E_2$ linked by a cyclic isogeny of degree $m$.
\end{theorem}

\paragraph{Proof.}
Denote by $T_m \subset X_0(N)\times X_0(N)$ the Hecke correspondence of level
$m$, and let $C_m \subset E\times E$ denote its image under $\pi\times\pi$. We
view $C_m$ as a symmetrical correspondence on $E$. 

Suppose that $C_m$ contains infinitely many torsion points of the abelian
variety $A=E\times E$. Then it follows from the Manin-Mumford Conjecture,
proved by Raynaud (see \cite{Raynaud} for the relevant case), that $C_m$ is the
translate by a torsion point of an abelian subvariety of $A$.  Now, the
one-dimensional abelian subvarieties of $A$ are of the form $\{0\}\times E$, $E
\times \{0\}$, or graphs of endomorphisms of $E$. But $C_m$ is symmetrical,
hence it is a translate of the graph of an automorphism of $E$, so $C_m$ is a
correspondence of degree one. This implies that $\deg(T_m)\leq\deg(\pi)$, and
the result follows, as $\deg(T_m)=\psi(m)\geq m+1$.\qed


\end{document}